\documentclass[12pt,draftcls,onecolumn]{IEEEtran}

\usepackage{pdfsync}
\usepackage{graphicx}
\usepackage{amsmath,amssymb}
\usepackage{epstopdf}
\usepackage{amsfonts}
\usepackage{euscript}
\usepackage{mathrsfs}
\usepackage{marvosym}
\usepackage[pdftex]{color}
\usepackage{pdfcolmk}
\usepackage{url}

\definecolor{Royalblue}{cmyk}{1,0.30,0.2,0.2}

\usepackage[dvips]{epsfig} \DeclareGraphicsExtensions{.ps,.eps}
\def\I{\mathbb{I}}

\def\R{\mathbb{R}}
\def\C{\mathbb{C}}

\def\im{{\rm im}\,}

\def\rank{{\rm rank}\,}

\newcommand{\tp}{^{\top}}
\newcommand{\mtp}{^{-\top}}
\newcommand{\inv}{^{-1}}
\newcommand{\script}[1]{\EuScript{#1}}

\newcommand{\diag}{{\rm diag}}
\newcommand{\beq}{\begin{equation}}
\newcommand{\eeq}{\end{equation}}
\newcommand{\bea}{\begin{eqnarray}}
\newcommand{\eea}{\end{eqnarray}}

\newcommand{\bsea}{\begin{subeqnarray}}
\newcommand{\esea}{\end{subeqnarray}}
\newcommand{\nn}{\nonumber}

\newcommand{\proof}{\noindent {\it Proof. }}

\newcommand{\qed}{\hfill $\Box$ \vskip 2ex}
\newcommand{\rst}{\setcounter{acount}{0}}
\newcommand{\stoa}{\setcounter{acount}{1}}
\newcommand{\stob}{\addtocounter{equation}{-1} \setcounter{acount}{2}}
\newcommand{\stoc}{\addtocounter{equation}{-1} \setcounter{acount}{3}}
\newcommand{\stod}{\addtocounter{equation}{-1} \setcounter{acount}{4}}

\newcommand{\bmat}{\left[ \begin{matrix}}
\newcommand{\emat}{\end{matrix} \right]}

\newcounter{acount}

\newtheorem{thm}{Theorem}[section]
\newtheorem{ass}{Assumption}[section]
\newtheorem{deff}{Definition}[section]
\newtheorem{lem}{Lemma}[section]
\newtheorem{cor}{Corollary}[section]
\newtheorem{prob}{Problem}[section]
\newtheorem{prop}{Proposition}[section]
\newtheorem{rem}{Remark}[section]
\newtheorem{chal}{Challenge}[section]
\newtheorem{misi}{Mission Impossible}[section]

\newtheorem{cjt}{Conjecture}[section]

\newcommand{\bthm}{\begin{thm}}
\newcommand{\bdeff}{\begin{deff}}
\newcommand{\bass}{\begin{ass}}
\newcommand{\blem}{\begin{lem}}
\newcommand{\bcor}{\begin{cor}}
\newcommand{\bprob}{\begin{prob}}
\newcommand{\bprop}{\begin{prop}}
\newcommand{\brem}{\begin{rem}}
\newcommand{\bchal}{\begin{chal}}
\newcommand{\bmisi}{\begin{misi}}

\newcommand{\ethm}{\end{thm}}
\newcommand{\edeff}{\end{deff}}
\newcommand{\eass}{\end{ass}}
\newcommand{\elem}{\end{lem}}
\newcommand{\ecor}{\end{cor}}
\newcommand{\eprob}{\end{prob}}
\newcommand{\eprop}{\end{prop}}
\newcommand{\erem}{\end{rem}}
\newcommand{\echal}{\end{chal}}
\newcommand{\emisi}{\end{misi}}

\newcommand{\pinv}{^+}

\title{\Large Representation and Factorization of Discrete-Time Rational All-Pass Functions}

\author{Augusto Ferrante and Giorgio  Picci,~\IEEEmembership{Life~Fellow,~IEEE}% <-this % stops a space
\thanks{A. Ferrante and 
 G. Picci are with   the Department of Information
Engineering, University of Padova, Padova, Italy; e-mail: {\tt\small augusto@dei.unipd.it} (A. Ferrante); {\tt\small picci@dei.unipd.it} (G. Picci).}
}

\begin{document}
 \maketitle

%$\aleph \mathfrak{A}\mathfrak{a}\mathscr{A}$ \EUR   \EURhv \EURcr   \EURtm 

\begin{abstract}
We obtain a general characterization of discrete-time all-pass rational matrix functions from state-space representations. It can be employed to address model reduction problems in the same vein of the theory developed by Glover in the continuous-time. Besides model reduction, this characterization is shown to be useful in a variety of contexts such as studying LMI's and Riccati equations and especially in the factorization of all-pass functions. The results  are obtained in the most general setting,    without introducing any {\em ad hoc} assumption.
\end{abstract}
\thispagestyle{empty}

\section{Introduction}\label{Introd}

In this paper we  provide a completely general  characterization and parameterization of  discrete-time all-pass matrix  functions and use this result  to describe in full generality the geometry of the solution set of  certain  LMI's  and of the associated Riccati equations. We also develop   a factorization theory and related state-space procedures for the factorization of all-pass functions. The characterization  of  discrete-time all-pass matrix  functions  presented in   the main theorem  of  the next section,  parallels the continuous-time   fundamental  result of   Glover's  \cite[Theorem 5.1]{Glover-84} in the most general setting,    without   introducing any {\em ad hoc}  assumption. A detailed rigorous proof of this result  seems to be presented  here for the first time after   past unfruitful attempts in the literature.  Commonly used facilitating assumptions in discrete-time, such as non-singularity of various  matrices (in particular of the $A$ matrix) and unmixing  are  avoided and only discussed as particular corollaries of more general statements.\\
A discrete-time  version of Glover's model reduction procedure seems to be worth  as discrete-time models are often the rule in applications. Its derivation however  is  not a simple transposition  of the arguments used in  continuous-time   as there are several differences which make the job technically much harder, see for example the attempts in  \cite{Gu-05}.   Some of these difficulties are well-known.  The use of the often advocated Cayley transform  requires for example  invertibility assumptions  which are not met in some applications  (see, e.g., the   comment in \cite[p. 1996]{Wimmer-06}). 
Moreover,  it seems to be an accepted point in the systems and control community that,    as stated in \cite[p. 559]{Zhou-etal}, {\em ``\dots it is generally more appealing to give derivations in the coordinates of the original [discrete-tme] data; also algorithms may be more reliable if generated for the specific model class''}. Apparently   a  discrete version of the continuous-time all-pass dilation of Glover under  general hypotheses as those made in the present  paper  has been lacking.   So far, to our best knowledge, the book literature of the last two or three decades, e.g. \cite{Alpay-G-88} or, \cite{Zhou-etal} \cite{Fuhrmann-book}  seems to be    just re-proposing   continuous-time $H^{\infty}$ model reduction   and  does not   address a discrete-time version of Glover's theory. 
     
The   results of the paper have many possible applications. Applications to Hankel-norm approximation of rational discrete-time transfer functions  may  now be pursued by just  following the route shown in the paper \cite{Glover-84}.   In Chapter 16 of the book \cite{LPBook}   a slightly less general characterization of discrete all-pass functions is used  to do   Hankel-norm stochastic model approximation. Stochastic modeling without stability constraints is another direction which has been touched upon in \cite {Ferrante-P-15}, further exposed in \cite{LPBook} and can be addressed   in wider generality by using the techniques described in this paper.   This is a relatively unappreciated area of stochastic modeling   which has several applications to smoothing and to non causal estimation. We believe that   this setting  is worth understanding especially because of a very illuminating  isomorphism with LQ control with an indefinite  cost function.  In a companion paper \cite{Ferrante-P-15} we shall    apply   this isomorphism  to resolve an old open problem about the existence of negative semidefinite solutions of the Riccati equation of LQ control. 

The lay-out of this paper is as follows:\\
Section \ref{MainSect} contains the statement and proof of the main result. The proof is essentially self-contained save for a technical Lemma from \cite{Ferrante-Ntog-Automatica-13} which considerably generalizes a result on controllability due to  Wimmer \cite{Wimmer-92}.
\\
 In section \ref{LMIRic} we introduce two dual linear matrix inequalities with a rank constraint which define families of square all-pass functions   having a fixed pole structure. We  prove a geometric characterization of all solutions in terms of $A$- or $A\tp$- invariant subspaces. When $A$ is non singular these matrix inequalities turn into two dual homogeneous algebraic Riccati equations. A very exhaustive classification and description of the solutions of those Riccati equations is provided. It is well-known, see e.g.  \cite{Wimmer-06} that  the analysis of algebraic Riccati equations can be   reduced to that of homogeneous Riccati equations.
\\
 The study of families of solutions of the constrained LMI's of Section  \ref{LMIRic} unveils the basic principles and a direct method to characterize and classify the left- and right all pass factors of an arbitrary square all pass rational function. Rational factorization theory was first systematically discussed in the early book \cite{BGK}  quite heavily relying on the assumption of an invertible $D$ matrix. Here we extend the factorization results of Fuhrmann and Hoffmann \cite{Fuhrmann-H-95} derived for inner functions, under general hypotheses. When $A$ is non-singular the classification can be given  directly in  terms of solutions of  two dual homogeneous algebraic Riccati equations.
\\
 In the concluding section  we indicate some possible generalizations to non square matrix functions.

Notations in the paper are quite    standard; we only mention that   $X\pinv$ denotes the 
  MooreÐPenrose pseudoinverse  of the matrix $X$.  A technical condition which is often referred to is that of {\em unmixing}. One says that $A\in  \R^{n\times n}$ has {\em unmixed spectrum} or, briefly, is {\em unmixed} if it does not have reciprocal pairs of eigenvalues. In particular an unmixed matrix cannot have  eigenvalues of modulus one.

\section{The main result}\label{MainSect}
\bthm \label{maintheorem-ap} 
\ 
\begin{enumerate}
\item\label{point1}
Let 
\beq\label{realization-of-Q}
 Q(z):= C(zI-A)^{-1}B+D
\eeq
be a minimal realization of an $m\times m$ rational discrete-time all-pass function. 
Then  $A$ is non-singular if and only if $D$ is non-singular.
\item\label{point2}
Let (\ref{realization-of-Q}) be a minimal realization of a  rational discrete-time all-pass function.
 Then there exist two invertible matrices $P=P\tp$ and $Q=Q\tp$ such that $PQ=I$ and
\beq\label{eqforP}
\left\{
\begin{array}{l}
A P A\tp - P =BB\tp\\
BD\tp-A P C\tp=0\\
DD\tp-CPC\tp=I
\end{array}
\right.
\eeq
\beq\label{eqforQ}
\left\{
\begin{array}{l}
A\tp Q A - Q =C\tp C\\
C\tp D-A\tp Q B=0\\
D\tp D-B\tp QB =I
\end{array}
\right.
\eeq
\item\label{point2.1}
Let (\ref{realization-of-Q})  be a minimal realization of a  rational discrete-time all-pass function. If equations  (\ref{eqforP}) admit a solution $P$, then such a $P$ is unique. If equations  (\ref{eqforQ}) admit a solution $Q$, then such a $Q$ is unique.
 
\item\label{point3}
Let $A\in\R^{n\times n},B\in\R^{n\times m},C\in\R^{m\times n},D\in\R^{m\times m}$ be given (no minimality is now assumed).
If there exists $P=P\tp$ satisfying (\ref{eqforP})  then $Q(z)$ given by (\ref{realization-of-Q}) is all-pass.\\
Similarly, if  there exists $Q=Q\tp$ satisfying (\ref{eqforQ})   then $Q(z)$ given by (\ref{realization-of-Q}) is all-pass.\\
 Finally, $P$ is a non-singular solution of (\ref{eqforP}) if and only if $P^{-1}$ is a (non-singular) solution of (\ref{eqforQ}).
\item\label{point4}
Let $A\in\R^{n\times n},B\in\R^{n\times m}$ be a given reachable pair. Then, $P=P\tp$ is such that
\beq\label{eqforPsolon}
A P A\tp - P =BB\tp
\eeq
if and only if there exist matrices $C\in\R^{m\times n}$ and $D\in\R^{m\times m}$  such that 
$Q(z)$ given by (\ref{realization-of-Q}) is a minimal realization of an all-pass function and 
$P$ is the solution of (\ref{eqforP}) for the quadruple $(A,B,C,D)$.
In this case, $P$ is necessarily non-singular and such that
\beq\label{propofsolpsdp}
I+B\tp P^{-1} B\geq 0.
\eeq
\item\label{point5}
Let $A\in\R^{n\times n},C\in\R^{m\times n}$ be a given observable pair.
Then,   $Q=Q\tp$ is such that
\beq\label{eqforQsolon}
A\tp Q A - Q =C\tp C
\eeq
if and only if there exist matrices $B\in\R^{n\times m}$ and $D\in\R^{m\times m}$  such that 
$Q(z)$ given by (\ref{realization-of-Q}) is a minimal realization of an all-pass function and  $Q$ is the solution of (\ref{eqforQ}) for the quadruple $(A,B,C,D)$.
In this case, $Q$ is necessarily non-singular and such that
\beq\label{propofsolqsdp}
I+C Q^{-1} C\tp \geq 0.
\eeq
\item \label{7}
Let $A\in\R^{n\times n},B\in\R^{n\times m},C\in\R^{m\times n}$ be given.
If there exists $P=P\tp$ and $Q=Q\tp$ such that
\beq\label{eqforPQ}
\left\{
\begin{array}{l}
A P A\tp - P =BB\tp\\
A\tp Q A - Q =C\tp C\\
PQ =I
\end{array}
\right.
\eeq
then there exists a matrix $D\in\R^{m\times m}$  such that 
$Q(z)$ given by (\ref{realization-of-Q}) is all-pass.
\end{enumerate}
\end{thm}

\proof

\noindent 1)\hspace{1.mm}
By assumption we have
\beq\label{app-n}
Q(z)Q^\ast(z)=I.
\eeq
Notice that $Q(\infty)=D$ so that by taking the limit  $z\rightarrow\infty$ in (\ref{app-n}), we see that $D$ is non-singular if and only if
$Q^\ast(z)$ is bounded at infinity or, equivalently, if and only if $Q(z)$ is bounded in a neighborhood of the origin. By taking into account that (\ref{realization-of-Q}) is a minimal realization, this is equivalent to $A$ being non-singular.

\noindent 2)\hspace{1.mm}
Let us first assume that $D$ is non-singular.
By recalling point \ref{point1}), we have that  $A$ is non-singular as well.

We have the following minimal realizations:
\beq
Q(z)^{-1}=D^{-1}-D^{-1}C (zI-\Gamma)^{-1} BD^{-1},\  \Gamma:= A-BD^{-1}C.
\eeq
and
\bea\nn
Q^\ast(z)&=&B\tp (z^{-1}I-A\tp)^{-1}C\tp +D\tp\\
&=&D_0\tp
-B\tp A^{-\top}(zI-A^{-\top})^{-1}A^{-\top}C\tp,
\eea
with $D_0\tp:=D\tp-B\tp A^{-\top}C\tp$,
so that, by   imposing $Q(z)^{-1}=Q^\ast(z)$, we conclude that there exists a unique invertible matrix $T$ such that
\bea\label{eqiuv}
\stoa\label{eqiuva}
&&
T^{-1}A^{-\top}T=A-BD^{-1}C (=\Gamma)\\
\stob\label{eqiuvb}
&&T^{-1}A^{-\top}C\tp=BD^{-1}\\
\stoc\label{eqiuvc}
&&B\tp A^{-\top} T=D^{-1}C\\
\stod\label{eqiuvd}
&&D^{-1}=D\tp-B\tp A^{-\top}C\tp
\eea\rst
By inserting (\ref{eqiuvc}) in (\ref{eqiuva}) and multiplying on the right side by $(A^{-\top} T)^{-1}$, we get
\beq\label{solprimaid}
T^{-1} -AT^{-1}A\tp =-BB\tp
\eeq
so that the first of (\ref{eqforP}) admits a solution $P=T^{-1}$.
Moreover, by inserting the expression of $D\inv$ provided by  (\ref{eqiuvd}) 
in (\ref{eqiuvb}) we get $B D\tp = (T\inv +BB\tp )A\mtp C\tp$, which, in view of 
(\ref{solprimaid}), may be written as 
$BD\tp=AT^{-1}C\tp$, so that $P=T^{-1}$ solves also the second of (\ref{eqforP}).
Finally, by multiplying  (\ref{eqiuvd}) on the left side by D and taking into account of
  (\ref{eqiuvc}), we easily see that $P=T^{-1}$ solves also the third of (\ref{eqforP}).
  Similarly we see that from (\ref{eqiuv}) it follows that $T$ solves the three equations
(\ref{eqforQ}).
The proof that $T$ is symmetric is a bit lengthy and is deferred to Appendix \ref{B}.

So far we have established our result in the case when $D$ is non-singular.
We now show how this case may be viewed as a first step for proving the result in the general setting in which $D$ may be singular.
Consider an arbitrary rational proper all pass function $Q(z)$ and the corresponding factorization (\ref{facaps+ns}) established in Lemma \ref{primolemmafactoriallp} of the appendix.
Let $Q_0(z):= C_0(zI-A_0)^{-1}B_0+D_0$ be a minimal realization of $Q_0(z)$ so that $D_0=Q_0(\infty)$ is non-singular.
Then equations (\ref{eqforP}) with $A=A_0$, $B=B_0$, $C=C_0$ and $D=D_0$ have a symmetric solution $P_0$ which is non-singular.
In view of Lemma \ref{secondolemfacapf} we know that $Q_1(z):=Q_0(z)\bar{Q}_1(z)$ has the  reachable realization $Q_{1}(z)=C_1\left(zI-A_1\right)^{-1}B_1 +D_1$ where 
$$
C_1:=[D_{0,2}\mid C_0],  A_1:=\bmat0&0\\B_{0,2}&A_0\emat, B_1:=\bmat0 & I\\B_{0,1}&0\emat U_{1},$$
and
$D_1:=[D_{0,1}\mid 0]U_{1}$.
Now it is immediate to check by inspection that 
\beq
P_1:=\bmat I&0\\0&P_0\emat
\eeq
solves equations (\ref{eqforP}) with $A=A_1$, $B=B_1$, $C=C_1$ and $D=D_1$.
We can iteratively repeat this argument for $Q_i(z)$, $i=2,3,\dots,k$ and  eventually find that 
$Q(z)$ has a reachable realization $Q(z)=\bar{C}\left(zI-\bar{A}\right)^{-1}\bar{B}+D$
and that equations (\ref{eqforP}) with $A=\bar{A}$, $B=\bar{B}$, $C=\bar{C}$ and $D=D$, have a solution $\bar{P}$.
Without loss of generality we may assume that $\bar{A}$, $\bar{B}$, $\bar{C}$ are in the Kalman reachability form 
\beq
\bar{C}=[\tilde{C}\mid 0],\qquad \bar{A}=\bmat \tilde{A}&0\\A_{21}&A_{22}\emat,\qquad \bar{B}:=\bmat \tilde{B}\\B_{2} \emat
\eeq
and $\bar{P}$ is partitioned conformably as 
\beq
\bar{P}=\bmat \tilde{P}&P_{12}\\P_{12}\tp &P_{22}\emat.
\eeq
By writing  block-wise equations (\ref{eqforP}) with $A=\bar{A}$, $B=\bar{B}$, $C=\bar{C}$, $D=D$,  and $P=\bar{P}$ we  see that the $(1,1)$ block $\tilde{P}$  is a symmetric solution of equations (\ref{eqforP}) with $A=\tilde{A}$, $B=\tilde{B}$, $C=\tilde{C}$, $D=D$ corresponding to the minimal realization 
\beq\label{minrealQtilde}
Q(z)=\tilde{C}\left(zI-\tilde{A}\right)^{-1}\tilde{B}+D.
\eeq
The original minimal realization of $Q(z)$ is related to (\ref{minrealQtilde}) by a change of basis so that there exists a non-singular matrix $T$ such that 
$A=T\inv \tilde{A} T$, $B=T\inv \tilde{B}$, $C=\tilde{C}T$.
Then it is immediate to check that $P:=T\inv \tilde{P}T\mtp$ is a solution of equations (\ref{eqforP}) for the original realization (\ref{realization-of-Q}).
Observe that by minimality of the realization $(A,B,C)$ we have that $P$, solving the Lyapunov equation in  (\ref{eqforP}), is  non-singular.

By resorting to a dual   argument we establish the existence of a non-singular matrix $Q=Q\tp$ solving (\ref{eqforQ}).

It remains to show that $PQ=I$.
To this aim,  write (\ref{eqforP}) in the form
\beq
FXF\tp=X
\eeq
where
\beq\label{deffx}
F:=
\bmat
A&B\\C&D\emat \qquad X:=\bmat P & 0\\0&-I\emat.
\eeq 
Clearly, $X$ is non-singular and 
\beq\label{formofXinv}
X\inv =\bmat P^{-1} & 0\\0&-I\emat.
\eeq
Thus
$FXF\tp X\inv=I$.
Therefore, 
 $F$ is non-singular as well and we have $F\tp X\inv=X\inv F\inv$ or
 \beq\label{puqinv}
 F\tp X\inv F=X\inv.
 \eeq
The expression (\ref{formofXinv}) of $X\inv$  implies that $P^{-1}$
is a solution  of equations (\ref{eqforQ}).
As proven below these equations  however  admit a unique solution so that $P^{-1}=Q$ or, equivalently, $PQ=I$.

\noindent 3)\hspace{1.mm}
Assume that $P_1$ and $P_2$ are solutions of (\ref{eqforP}) and let $\Delta:=P_1-P_2$.
We need to show that $\Delta=0$.
It is immediate to check that $\Delta$ satisfies the equations
\beq\label{eqforDelta}
\left\{
\begin{array}{l}
A \Delta A\tp - \Delta =0\\
A\Delta  C\tp=0\\
C\Delta C\tp=0
\end{array}
\right.
\eeq
From the second and the third of these equations we  see that $\im(\Delta C\tp )$ is contained in the non-observability subspace of $(A,C)$ and, since $(A,C)$ is assumed to be observable, this means that 
$\Delta C\tp =0$.
This implies that $C\Delta =0$ and, in turn, $CA^k \Delta (A\tp)^k=0$ for all $k=0,1,\dots,n-1$, so that $CA^k \Delta (A\tp)^{n}=0$ for all $k=0,1,\dots,n-1$.
This means that $\im(\Delta (A\tp)^{n})$ is contained in the non-observability subspace of $(A,C)$, so that, as before, $\Delta (A\tp)^{n}=0$. Now, by multiplying the first of (\ref{eqforDelta}) on the right side by $(A\tp)^{n-1}$, we get $\Delta (A\tp)^{n-1}=0$ and, iteratively, $\Delta (A\tp)^{n-2}=0$, and so on, up to $\Delta =0$.
The proof for equations (\ref{eqforQ}) is dual and is therefore skipped.

\noindent 4)\hspace{1.mm}
Assume that equations (\ref{eqforP}) admit a solution $P=P\tp$.
Let us  compute the product
\bea
\nn
\Phi\!\!\!&:=\!\!&Q(z)Q\tp(z^{-1})\\
\nn
&=\!\!&[C(zI-A)^{-1}B+D][B\tp(z^{-1}I-A\tp)^{-1}C\tp+D\tp].
\eea
The first of equations (\ref{eqforP}) can be rewritten as
\bea\nn
BB\tp&=&(zI-A)P(z^{-1}I-A\tp)-zP(z^{-1}I-A\tp) \\
\nn
&&- z^{-1}(zI-A)P,
\eea
so that  
$$C(zI-A)^{-1}BB\tp(z^{-1}I-A\tp)^{-1}C\tp=$$
$$CPC\tp -zC (zI-A)^{-1} P C\tp-
z^{-1}CP(z^{-1}I-A\tp)^{-1} C\tp.$$
Moreover, from
$$z(zI-A)^{-1}=I+A(zI-A)^{-1} =I+(zI-A)^{-1}A $$
it follows  that
$$C(zI-A)^{-1}BB\tp(z^{-1}I-A\tp)^{-1}C\tp=$$
$$-CPC\tp -C(zI-A)^{-1}A P C\tp-
C P A\tp (z^{-1}I-A\tp)^{-1}C\tp.$$

In conclusion, we have
\bea
\nn
Q(z)Q\tp(z^{-1})&=& DD\tp-CPC\tp \\
\nn&&
+C(zI-A)^{-1}(BD\tp-A P C\tp)\\
\nn&&
+
(DB\tp-C P A\tp) (z^{-1}I-A\tp)^{-1}C\tp.
\eea
By taking into account the second and the third of equations (\ref{eqforP}),
we now get  $Q(z)Q\tp(z^{-1})=I$.

Assume now that equations (\ref{eqforQ}) admit a solution $Q=Q\tp$.
By computing the product $Q\tp(z^{-1})Q(z)$ and using the dual of the previous argument, we  get
$Q\tp(z^{-1})Q(z)=I$.

The fact that $P$ is a non-singular solution of (\ref{eqforP}) if and only if $P^{-1}$ is a non-singular solution of (\ref{eqforQ}) can be shown by defining $F$ and $X$ as in (\ref{deffx}) and using the same argument that led to (\ref{puqinv}).

\noindent 5)\hspace{1.mm}
One direction is an immediate consequence of point \ref{point2}).
For the converse,
since $(A,B)$ is, by assumption,  reachable, the solution $P$ of (\ref{eqforPsolon}) is invertible
\cite[Lemma 3.1]{Ferrante-Ntog-Automatica-13}.
Let $(n_+,n-n_+,0)$ be the inertia of $P$ which is equal to the inertia of $Q:=P^{-1}$.
Let
$
E:= \bmat -Q & 0 & A\tp\\0 & I_m & B\tp\\ A & B & -Q^{-1}\emat.$
The inertia of $E$ is 
given by the inertia of
$ \bmat -Q & 0 \\0 & I_m \emat$, i.e. $(m+n-n_+,n_+,0)$, plus the inertia of the corresponding Schur complement $S$ which is given by
\bea
\nn
S&:=&-Q^{-1}-[A\ \ B]\bmat -Q & 0 \\0 & I_m \emat^{-1}\bmat A\tp \\B\tp \emat\\
\nn
&=&-P+A P  A\tp -BB\tp=0_{n\times n}.
\eea
In conclusion, the inertia of $E$ is  $(m+n-n_+,n_+,n)$.
On the other hand, the inertia of $E$ is also given by 
the inertia of $-Q^{-1}=-P$, i.e. $(n-n_+,n_+,0)$, plus the inertia of the corresponding Schur complement $W$ which is given by
\bea\nn
W&:=&\bmat -Q & 0 \\0 & I_m \emat-\bmat A\tp \\B\tp \emat (-Q^{-1})^{-1}[A\ \ B]\\
&=&
\label{defWpsd}
\bmat
A\tp Q A - Q & A\tp Q B\\
B\tp QA & B\tp QB +I\emat.
\eea
Hence the inertia of $W$ is given by the inertia of $E$, i.e. $(m+n-n_+,n_+,n)$ minus the inertia of $-Q^{-1}=-P$, i.e. $(n-n_+,n_+,0)$, which amounts to $(m,0,n)$.
Thus, $W$ is positive semidefinite and has rank equal to $m$.
Therefore, there exists a full row-rank matrix
$
[C \mid D]\in\R^{m\times (n+m)}
$ 
such that
$
W =[C \mid D]\tp [C \mid D].$
This means that for the given  $A$ and $B$ and for the $C$ and $D$ obtained by previous developments,
there exists a $Q=P^{-1}$ solving (\ref{eqforQ}).
Then, in view of point \ref{point3}), the corresponding $Q(z)$ given by (\ref{realization-of-Q}) is all-pass.

We now prove (\ref{propofsolpsdp}). Indeed, we have already proved that $P^{-1}$ solves (\ref{eqforQ}) and from the third of these equations (\ref{propofsolpsdp}) follows immediately.

It remains to show that $(A,C)$ is an observable  pair. To address this issue we exploit (\ref{eqforQ}) whose validity we have already proven.
Assume now by contradiction that the pair $(A,C)$ is not observable and
let $V$ be a  full column-rank matrix whose columns (at least one by the contradiction assumption) form a basis for the unobservable subspace $\mathscr{N}:=\ker \bmat C\\CA\\ \vdots\\
CA^{n-1}\emat$, so that 
\beq
CV=0
\eeq
 and since $\mathscr{N}$ is $A$-invariant, there exists a matrix $K$ such that 
 \beq
 AV=VK.
 \eeq
By multiplying the first of  (\ref{eqforQ}) on the right side by $V$ we get $A\tp QAV=QV$.
We now multiply the first of  (\ref{eqforQ}) on the right side by $AV$ and on the left side by $A\tp$: We get $(A\tp)^2 QA^2V=A\tp QAV=QV$.
We can iterate this argument and multiply the first of  (\ref{eqforQ}) on the right side by $A^kV$ and on the left side by $(A\tp)^{k}$, $k=2,3,\dots$,   getting
\beq\label{iterazionediLyap}
(A\tp)^l QA^lV=QV,\qquad l=1,2,\dots.
\eeq
We now show that 
\beq\label{defdiU}
U:= Q A^n V\neq 0,
\eeq
where $n$ is the dimension of $A$. In fact, from (\ref{iterazionediLyap}) we get $(A\tp)^n U=(A\tp)^n QA^nV=QV$  and since $Q$ is non-singular and $V$ has full column-rank this yields (\ref{defdiU}).
We now consider the second of equations (\ref{eqforQ}). From this  equation, we get
$
D\tp C= B\tp QA,
$
and by right-multiplication by $V$, we get 
\beq
B\tp QAV=0
\eeq
so that 
\beq
B\tp QA^l V=B\tp QA V K^{l-1} =0,\qquad l=1,2,\dots, n-1.
\eeq 
Thus, for any $l=0,1,\dots, n-1$, we have
\bea
\nn
B\tp (A\tp)^l U &\!\!\!\!=\!\!\!\!&  B\tp (A\tp)^l Q A^n V= B\tp (A\tp)^l Q A^l V K^{n-l}\\
&\!\!\!\!=\!\!\!\!&  B\tp   Q   V K^{n-l}= B\tp   Q   A V K^{n-l-1}=0.
\eea
In conclusion, $\im(U)\neq\{0\}$ is contained in the unobservable subspace of the pair $(A\tp, B\tp)$ and this is a contradiction because $(A,B)$ is, by assumption, reachable, so that $(A\tp, B\tp)$ is observable.

\noindent 6)\hspace{1.mm}
This point is the dual of the previous one.

\noindent 7)\hspace{1.mm}
Since $P$ is clearly invertible we can use the same argument employed in the proof of 
point \ref{point4}) to show that
\bea\nn
W&:=&\bmat -Q & 0 \\0 & I_m \emat-\bmat A\tp \\B\tp \emat (-Q^{-1})^{-1}[A\ \ B]\\
\nn
&=&
\bmat
A\tp Q A - Q & A\tp Q B\\
B\tp QA & B\tp QB +I\emat
\eea
 is positive semidefinite and has rank equal to $m$.
Therefore, there exists a full row-rank matrix
$[C_0\mid D_0]\in\R^{m\times (n+m)}
$ such that
$
W =[C_0\mid D_0]\tp [C_0\mid D_0].$ 
In particular,
$$
A\tp Q A - Q =C\tp C=C_0\tp C_0
$$ 
so that there exists an orthogonal matrix $U$ such that $C=UC_0$.
Let $D:=UD_0$. Therefore,
\bea\nn
W&:=& [C_0\mid D_0]\tp [C_0\mid D_0]
\\
\nn
&=&
[C_0\mid D_0]\tp U\tp U [C_0\mid D_0]
=[C\mid D]\tp   [C\mid D].
\eea
In conclusion, we have
\beq
D\tp D=I+B\tp Q B
\eeq
and \beq
D\tp  C=B\tp QA
\eeq
These two equations together with the second of (\ref{eqforPQ}) give
(\ref{eqforQ}) and hence, in view of point \ref{point3}), 
$Q(z)= C(zI-A)^{-1}B+D$ is all-pass.
\qed

\brem
  In point \ref{point4}) of Theorem \ref{maintheorem-ap} the assumption of reachability of $(A,B)$   can probably be  eliminated for the first part of the result.
More precisely, we suggest the following   conjecture whose proof, however, seems to be non-trivial. \erem
\begin{cjt} 
Let $A\in\R^{n\times n},B\in\R^{n\times m}$ be given.
Then, there exists $P=P\tp$ such that
$
A P A\tp - P =BB\tp
$
if and only if there exist matrices $C\in\R^{m\times n}$ and $D\in\R^{m\times m}$  such that 
$Q(z)$ given by (\ref{realization-of-Q}) is an observable realization of an all-pass function.
\end{cjt}
Of course, we have a dual conjecture for point \ref{point5}).

\brem
 Consider an all-pass function $Q(z)$ represented by (\ref{realization-of-Q}).
Clearly $Q(z) U$ is still all-pass for any orthogonal matrix $U$.
The two functions $Q(z)$ and $Q(z) U=C(zI-A)^{-1}BU+DU$ correspond to the same dynamics so that it is natural to regard these two functions as equivalent. By considering the polar decomposition of $D$ we immediately see that for any given $D$ there is a unique $D_0=DU$ such that $D_0=D_0\tp\geq 0$.
Therefore, from now on, whenever convenient, we can safely assume, without loss of generality, that the ``$D$" matrix of the all-pass function $Q(z)$ is symmetric and positive semidefinite.
\erem

\brem
Consider point \ref{point4})  (or \ref{point5})) of Theorem \ref{maintheorem-ap}.
If $A$ is unmixed, once given $A$ and $B$, the solution $P$ of (\ref{eqforPsolon}) is uniquely determined and hence 
also the matrices $C$ and $D$ for which $Q(z) =C(zI-A)^{-1}B+D$ is all-pass are uniquely determined up to multiplication on the left side by a common orthogonal matrix. This is not the case when $A$ is not unmixed. 
In this case, for any particular solution $P$ of (\ref{eqforPsolon}) there exists a particular pair of matrices $C$ and $D$ (essentially different, i.e. not differing for multiplication on the left side by a common orthogonal matrix) for which $Q(z) =C(zI-A)^{-1}B+D$ is all-pass.
Notice, however, that, once fixed $A$, $B$ {\em and} $P$, the matrices $C$ and $D$  are always uniquely determined up to multiplication on the left side by a common orthogonal matrix.\\ 
Similar considerations can be made for \ref{point5}).
For example, let 
$A=\bmat 2 &0\\0&1/2 \emat$ and $C=I$.
In this case the set of all solutions $Q$  of (\ref{eqforQsolon}) can be parametrized as $Q=\bmat 1/3 &q\\q& -4/3\emat$ with $q$ being a real parameter.  For example, for $q=0$, we get
$B_0=\bmat 3 &0\\0& -3/4\emat$ and $D_0=\bmat 2 &0\\0& 1/2\emat $, where the degree of freedom corresponding to the choice of an arbitrary orthogonal matrix multiplying both $B_0 $ and $D_0 $ on the right side, has been fixed in such a way that $D_0=D_0\tp\geq 0$: since $D_0$ is non-singular this procedure does not leave any further degree of freedom.
For $q=1/6$, we get $B_{1/6}=\bmat 2.85 &0.57\\0.14& -0.71\emat$ and $D_{1/6}=\bmat 1.95 &0.14\\0.14& 0.52\emat$, where, again, the degree of freedom corresponding to the  arbitrary orthogonal matrix  has been fixed in such a way that $D_{1/6}=D_{1/6}\tp\geq 0$.
In conclusion, the two solutions corresponding to $q=0$ and $q=1/6$ lead to all-pass functions with different dynamical properties.
\erem

\section{LMI's and Homogeneous Algebraic Riccati equations}\label{LMIRic}
All-pass functions can be seen as  spectral factors of a spectral density function identically equal to the identity matrix; i.e. $\Phi(z)\equiv I$. This point of view turns out to be  useful for classification of all-pass functions having a pre-assigned pole dynamics. It is a classical result in system and control theory \cite{Willems-71} that rational spectral factorization can be cast in terms of linear matrix inequalities (LMI). This point of view will be used here. It will be further developed in a forthcoming  companion paper \cite{Part2} devoted to stochastic modeling.  In this section we shall just consider square spectral factors which are all-pass.

To fix the pole dynamics we may either assign a reachable pair $(  A\in\R^{n\times n},\, B\in\R^{n\times m})$ or an observable pair $(\, C\in\R^{m\times n},\, A\in\R^{n\times n})$. These are two ``dual'' structural data which will be fixed hereafter. Accordingly, define
$$
M(P):=\bmat
A P A\tp - P & A P C\tp \\
C PA\tp & C P C\tp +I\emat,
$$
and
$$
N(Q):=\bmat
A\tp Q A - Q & A\tp Q B\\
B\tp QA & B\tp QB +I\emat
$$
and consider the two  dual, constrained linear matrix inequalities (CLMI),
\beq
\label{garep}
\left\{\begin{array}{l}
M(P)\geq 0\\
{\rm rank}[M(P)]=m
\end{array}
\right.
\eeq
\beq
\label{gareq}
\left\{\begin{array}{l}
N(Q)\geq 0\\
{\rm rank}[N(Q)]=m
\end{array}
\right.
\eeq
The following is an immediate corollary of Theorem \ref{maintheorem-ap}.
\begin{cor}
Let $P=P\tp$ be a solution of \eqref{garep} and let
\beq\label{eqforPLMI}
M(P) =\bmat G\\L\emat \bmat G\tp & L\tp \emat
\eeq 
be  a factorization of  full rank $m$. Then 
\beq
Q_{L}(z):= C(zI-A)^{-1}G +L
\eeq
is a (in general non minimal) realization of a square all-pass function. Dually, let
 $Q=Q\tp$ be a solution of \eqref{gareq} and let
\beq\label{eqforQLMI}
N(Q)=\bmat H & J\emat \tp \bmat H  & J\emat
\eeq 
be  a factorization of  full rank $m$. 
Then 
\beq
Q_{R}(z):= H(zI-A)^{-1}B +J
\eeq
is a  realization of a square  all-pass function.
\end{cor}

\smallskip

Clearly $P=0$ and $Q=0$ are always solutions of the inequalities \eqref{garep} and \eqref{gareq} and it may well happen  that these inequalities admit no other solutions save for these trivial ones. We need to exclude these uninteresting circumstances. We shall henceforth assume that there is a $D$ such that the matrix function with (minimal) realization
\beq \label {realization-of-Q-2new}
Q(z):=  C(zI-A)^{-1}B +D
\eeq
is all-pass. By Theorem \ref{maintheorem-ap} this happens  if and only  equations \eqref{eqforPLMI} with $G=B$ and $L=D$ hold for a nonsingular $P\equiv P_0$ or, equivalently, if and only if \eqref{eqforQLMI} with   $H=C$ and $J=D$ hold for a nonsingular $Q\equiv Q_0$.  $P_0$ and $Q_0$ turn in fact out to be  such that $P_0Q_0=I$.
In the next section it will be shown that each $Q_R(z)$ is a right factor of $Q(z)$ and each $Q_L(z)$ is a left factor of $Q(z)$.\\
{\em Notational convention:}   From now on, $(A,B,C,D)$ such that (\ref{realization-of-Q-2new}) is a minimal realization of a square all-pass function will be the problem data;  the unique solutions of \eqref{eqforP} and  \eqref{eqforQ} will be denoted by $P_0$ and $Q_0$, respectively and we shall reserve the symbols $P$ and $Q$ for generic solutions  of  \eqref{garep}  and \eqref{gareq}.\\ 

\bthm\label{thmhgare}
Let (\ref{realization-of-Q-2new}) be a minimal realization of a square all-pass function. Then
\begin{enumerate}
\item
(i) For each solution $P=P\tp$ of  (\ref{garep}),  the    subspace
\beq
{\script Y}=\ker(P)
\eeq
is $A\tp$-invariant.\\
(ii) 
The set  of  non-singular solutions of (\ref{garep}) can be parametrized as: 
\beq
{\mathbb P}= \{ P_{\Delta}: \ \Delta\in{\mathcal D}_p\}\,,
\eeq
where $P_{\Delta}:= (P_0^{-1}+\Delta)^{-1}$, $P_0$ is the unique solution of (\ref{eqforP}), and  ${\mathcal D}_p$ is the vector space  of solutions of $A\tp \Delta A -\Delta=0$.
If $A$ is unmixed, then ${\mathcal D}_p=\{0\}$ and (\ref{garep}) admits a unique non-singular solution $P_{\Delta}=P_0$, which is the unique   solution  of (\ref{eqforP}).
If $A$ is not unmixed, then ${\mathbb P}$ contains infinitely many solutions.
\\
(iii)  Let $ P_{\Delta}$ be a non-singular solution of (\ref{garep});   then to any $A\tp$-invariant subspace ${\script Y}$   there corresponds a solution $P$ of (\ref{garep}) given by 
\beq\label{paramgenforsetutte}
P:=\left[(I-\Pi)P_{\Delta}^{-1} (I-\Pi)\right]\pinv
\eeq
where $\Pi$ is the orthogonal projector onto ${\script Y}$. The kernel of $P$ is
${\script Y}$. 
If $A$ is unmixed, equation (\ref{paramgenforsetutte}), with $P_{\Delta}=P_0$ being the unique solution of (\ref{eqforP}), parametrizes the set of all   solutions of (\ref{garep}) in terms of $A\tp$-invariant subspaces.

\item
(i) For each solution $Q=Q\tp$ of  (\ref{gareq}),  the subspace  
\beq
{\script X}=\ker(Q).
\eeq
is   $A$-invariant.\\
(ii) 
The set ${\mathbb Q}$ of non-singular solutions of (\ref{gareq}) can be parametrized as: 
\beq\label{tuttelesolinvertibili}
{\mathbb Q}= \{Q_{\Delta}: \ \Delta\in{\mathcal D}_q \}
\eeq
where $Q_{\Delta}:=(Q_0^{-1}+\Delta)^{-1}$, $Q_0$ is the unique solution of (\ref{eqforQ}), and ${\mathcal D}_q$ is the vector space  of solutions of $A\Delta A\tp -\Delta=0$.
If $A$ is unmixed, then ${\mathcal D}_q=\{0\}$ and (\ref{gareq}) admits a unique non-singular solution $Q_{\Delta}=Q_0$, which is the unique solution    of (\ref{eqforQ}).
If $A$ is not unmixed, than ${\mathbb Q}$ contains infinitely many solutions.
\\
(iii)    Let  $Q_{\Delta}$ be a non-singular solution of (\ref{gareq}),  then to any  $A$-invariant subspace ${\script X}$, there corresponds a solution $Q$ of (\ref{garep}) given by 
\beq\label{paramgenforsetutteq}
Q:=\left[(I-\Pi) Q_{\Delta}^{-1} (I-\Pi)\right]\pinv
\eeq
where $\Pi$ is the orthogonal projector onto ${\script X}$.
The kernel of $  Q$ is ${\script X}$. 
If $A$ is unmixed, equation (\ref{paramgenforsetutteq}), with $Q_{\Delta}=Q_0 $ being the unique solution of (\ref{eqforQ}), parametrizes the set of all  solutions of (\ref{gareq}) in terms of $A$-invariant subspaces.
\end{enumerate}
\ethm
\proof
We prove only point 2), as the proof of point 1) is dual.\\
(i) It is clear that (\ref{gareq}) is equivalent to existence of two matrices $H\in \R^{m\times n}$ and $J\in\R^{m\times m}$ such that $[H \mid J]$ has full row-rank and $N(Q)=[H \mid J]\tp[H \mid J]$.
Therefore, if $Q$ is a solution of (\ref{gareq}) then $A\tp Q A-Q=H\tp H $, so that, in view of \cite[Lemma 3.1]{Ferrante-Ntog-Automatica-13}, ${\script X}:=\ker(Q)$ is $A$-invariant. \\
(ii) Clearly the solution $Q_0$ of (\ref{eqforQ}) is a non-singular solution of (\ref{gareq}) and the corresponding matrices $H$ and $J$, introduced before, coincide with 
  $C$ and $D$ of \eqref{realization-of-Q-2new}. Then in view of Theorem \ref{maintheorem-ap}, point \ref{point3}), we have
\beq\label{lyapinversa}
AQ_0^{-1}A\tp -Q_0^{-1}=BB\tp.
\eeq
 Let now $\tilde{Q}_0$ be another non-singular solution of (\ref{gareq}) and $C_0$ and $D_0$ be such that $N(\tilde{Q}_0)=[C_0 \mid D_0]\tp[C_0 \mid D_0]$. Equivalently,  $\tilde{Q}_0$ is a non-singular solution of (\ref{eqforQ}) corresponding to the quadruple $(A,B,C_0,D_0)$. Using again Theorem \ref{maintheorem-ap}, point \ref{point3}), we have that $\tilde{Q}_0^{-1}$ is a solution of (\ref{eqforP}) corresponding to the same quadruple, so that, in particular,
$A\tilde{Q}_0^{-1}A\tp -\tilde{Q}_0^{-1}=BB\tp$.
Comparing the latter with (\ref{lyapinversa}), we see that $\tilde{Q}_0^{-1}=Q_0^{-1}+\Delta$ where $\Delta$ is a solution of the homogeneous Lyapunov equation $A\Delta A\tp -\Delta=0$.
If $A$ is unmixed, this  equation has a unique solution $\Delta=0$ so that $\tilde{Q}_0=Q_0$.\\
Assume now that $A$ is not unmixed. Then equation $A\Delta A\tp -\Delta=0$ has a non-trivial vector space ${\mathcal D}_q$ of solutions and the previous argument shows that any non-singular solution $Q_{\Delta}$  of (\ref{gareq}) has the form $(Q_0^{-1}+\Delta)^{-1}$. It remains to show that all the elements of ${\mathbb Q}$ are solutions   of (\ref{gareq}). Observe that $A[Q_0^{-1}+\Delta]A\tp -[Q_0^{-1}+\Delta]=BB\tp$ for any $\Delta\in{\mathcal D}_q$.
Since $(A,B)$ is reachable, any $Q_{\Delta}:=Q_0^{-1}+\Delta$ is invertible and, in view of 
Theorem \ref{maintheorem-ap}, point \ref{point4}), there exist two matrices $C_\Delta$ and  $D_\Delta$
such that $C_{\Delta}(zI-A)^{-1}B+D_{\Delta}$ is a minimal realization of a  rational all-pass function and therefore $P_{\Delta}= Q_{\Delta}^{-1}$ is the solution of (\ref{eqforP}) corresponding to  the quadruple $(A,B,C_{\Delta},D_{\Delta})$. This is equivalent to $Q_{\Delta}:=P_{\Delta}^{-1}=(Q_0^{-1}+\Delta)^{-1}$,
being the solution of (\ref{eqforQ}) for the same quadruple so that $Q_{\Delta}$ is a solution of (\ref{gareq}) which therefore has  infinitely many solutions.\\
(iii) 
Let ${\script X}$ be an $A$-invariant subspace.
Consider  an orthogonal  change of basis induced by the matrix $T=[V_{\perp}\mid V]$, where the columns of $V$ form a basis for ${\script X}$ and the columns of $V_{\perp}$ form a basis for ${\script X}^\perp$.
In this basis we have
\beq\label{s1newbasisnew}
T\tp {\script X}=\im\bmat 0\\I\emat 
\eeq
and 
 \beq\label{a-inv-sp-formnew}
\bar{A}:=T\inv A T=T\tp AT=\bmat A_1&0\\A_{21} &A_2\emat.
\eeq 
 Partition $\bar{B}:=T\inv B=T\tp B$ conformably as $\bar{B}=\bmat B_1\\B_2\emat$.
Let $Q_\Delta$ be a non-singular solution of (\ref{gareq}) and let  $C_\Delta$ and $D_\Delta$ be such that $N(Q_\Delta)=[C_\Delta \mid D_\Delta]\tp[C_\Delta \mid D_\Delta]$. Equivalently,   $Q_\Delta$ is the non-singular solution of (\ref{eqforQ}) corresponding to an all-pass function described by the  quadruple $(A,B,C_\Delta,D_\Delta)$. Hence, in the new 
basis  $\bar{Q}_\Delta:=T\tp Q_\Delta T $ is a non-singular solution of (\ref{eqforQ}) corresponding to the quadruple $(\bar{A},\bar{B},\bar{C}_\Delta,D_\Delta)$, with $\bar{C}_\Delta:=C_\Delta T$.
In view of Theorem \ref{maintheorem-ap}, point \ref{point3}), $\bar{Q}_\Delta^{-1}$ is a non-singular solution of (\ref{eqforP}) corresponding to the same quadruple. Partition such a $\bar{Q}_\Delta^{-1}$ conformably with $\bar{A}$ as
 \beq\label{qinvnewbasis-new}
  \bar{Q}_\Delta^{-1}=\bmat P_{11}&P_{12}\\P_{12}\tp &P_{22}\emat
 \eeq
and note that  it must in particular satisfy the first  equation of (\ref{eqforP}) so that the block of index (1,1)   must satisfy  the reduced Lyapunov equation
\beq\label{eqLyinv-rednew}
A_1 P_{11} A_1\tp=P_{11} +B_1B_1\tp.
\eeq
Since  the pair $(A,B)$ is reachable, the pair $(A_1,B_1)$ is reachable as well, so that from Theorem \ref{maintheorem-ap}, point \ref{point4}), it follows that $P_{11}$ is invertible and there exist $C_1$ and $D_1$ such that $P_{11}$ is the unique solution of (\ref{eqforP}) corresponding to a reduced quadruple $(A_1,B_1,C_1,D_1)$ and  hence, 
$P_{11}^{-1}$ is the unique solution of (\ref{eqforQ}) corresponding to the same quadruple. It is now a matter of direct computation to check that 
\beq\label{q01}
\bar{Q} :=\bmat P_{11}^{-1}&0\\0 &0\emat
 \eeq
is a solution of (\ref{eqforQ}) corresponding to the quadruple $(\bar{A},\bar{B},[C_1\mid 0],D_1)$. Therefore,  $Q :=T \bar{Q}  T\tp$ is a solution of (\ref{eqforQ}) corresponding to the quadruple $(A,B,[C_1\mid 0]T\tp,D_1)$ and  hence, it is also a solution of (\ref{gareq}).
The fact that $\ker[Q ]={\script X} $ is  direct consequence of (\ref{q01}). We need to show that (\ref{paramgenforsetutteq}) is a coordinate-free representation of
$Q $. 
By observing that  $T\tp\Pi T=\bmat 0&0 \\0  &I\emat$ 
and
\bea\nn
(I-\Pi) Q_\Delta^{-1} (I-\Pi)&\!\!\!\!=\!\!\!\!&(I-\Pi)TT\tp Q_\Delta^{-1}T T\tp (I-\Pi)\\
\nn
&\!\!\!\!=\!\!\!\!&(I-\Pi)T\bmat P_{11} &P_{12}\\P_{12}\tp &P_{22}\emat T\tp (I-\Pi),
\eea
it is a straightforward computation to show that
\beq
\left[(I-\Pi) Q_\Delta^{-1} (I-\Pi)\right]\pinv=T \bmat P_{11}^{-1}&0\\0&0\emat T\tp=Q .
\eeq
The last thing that remains to be proven is the fact that when $A$ is unmixed,
 {\em all}    solutions of (\ref{gareq}) are parametrized in terms of $A$-invariant subspaces by (\ref{paramgenforsetutteq}), with $Q_{\Delta}=Q_0$.
We have already  shown that in this case  (\ref{gareq}) has a unique non-singular solution which coincides with the unique solution $Q_0$ of (\ref{eqforQ}).
The representation of   the other (singular) solutions can be obtained by a  procedure similar to the one introduced  above. % the only difference being that it requires a preliminary reduction that eliminate the singular components. %[NOT CLEAR]
Indeed, assume that  $Q$ is a singular solution of (\ref{gareq}) and let $H$, and $J$ be such that $N(Q)=[H \mid J]\tp[H \mid J]$.
As already proved, $\ker [Q]$ is $A$-invariant so that we can perform a change of coordinates such that in the new basis $Q$ has the structure of the right-hand side of (\ref{q01}), with $P_{11}$ being a non-singular matrix, $A$ has the structure of the right-hand side of (\ref{a-inv-sp-formnew}), and $B=[B_1\tp\mid B_2\tp]\tp$ and $H=[H_1\mid H_2]$ are partitioned conformably.
It is now a matter of direct computation to check that $P_{11}^{-1}$ is a solution of   (\ref{eqforQ}) corresponding to the quadruple $(A_1,B_1,H_1,J)$ so that
$P_{11}$ is a solution of (\ref{eqforP}) corresponding to the same quadruple.
Hence, $P_{11}$ satisfies the Lyapunov equation $A_1 P_{11}A_1\tp -P_{11}=B_1 B_1\tp$.
But since $A$ is unmixed, $A_1$ is also unmixed so that $P_{11}$ is uniquely determined by the Lyapunov equation. As a consequence, there is a unique $Q $ with the given kernel which necessarily coincides with the one given by the right-hand side of (\ref{paramgenforsetutteq}) with $Q_{\Delta}=Q_0$ and ${\script X} =\ker[Q].$
\qed

\brem\label{remdefcalpq}
Let ${\mathcal P}_{\Delta}$ and ${\mathcal Q}_{\Delta}$ denote the set of solutions of \eqref{garep} and \eqref{gareq} described by \eqref{paramgenforsetutte} and \eqref{paramgenforsetutteq} for a specific $\Delta$. While   when $A$ is unmixed (and hence ${\mathcal D}_p={\mathcal D}_q=\{0\}$ so that we necessarily have $\Delta=0$) the families ${\mathcal P}_0$ and ${\mathcal Q}_0$ constitute  the entire set of solutions of the LMI's  \eqref{garep} and (\ref{gareq}), it is not clear if this also holds for the case of a mixed $A$ even if one takes the union of the sets  ${\mathcal P}_{\Delta}$ with respect to  $\Delta\in{\mathcal D}_p$ or the union of the sets ${\mathcal Q}_{\Delta}$ with respect to  $\Delta\in{\mathcal D}_q$. %we have the following
The  theorem provides a bijective correspondence between the family ${\mathcal Q}_0$ of the solutions of (\ref{gareq}) and the family of $A$-invariant subspaces.  When $A$ is not unmixed, (\ref{gareq}), besides ${\mathcal Q}_0$, has infinitely many other families of solutions each of  which being likewise  parametrized by $A$-invariant subspaces. Each of these families corresponds to a non-singular solution $Q_{\Delta}\in{\mathbb Q}$ of (\ref{gareq}) where ${\mathbb Q}$ is the set of non-singular solutions parametrized by (\ref{tuttelesolinvertibili}).
The family ${\mathcal Q}_0$  corresponding to $\Delta=0$  will play  an important role  in the following.\\
Similar considerations can be made for the dual   family ${\mathcal P}_0$ of solutions of (\ref{garep}) which, in  case of unmixed $A$ constitutes the set of all solutions of (\ref{garep}) and in  case of a mixed $A$ is just one of   infinitely many families of solutions  of (\ref{garep}). 
\erem

\brem\label{rembijcpq}
There is an obvious  bijective correspondence between the set of $A$-invariant subspaces and that of $A\tp$-invariant subspaces. Indeed,   ${\script  X} $ is $A$-invariant if and only if  ${\script Y} :={\script  X}^\perp$ is  $A\tp$-invariant. This correspondence induces a bijective correspondence between the sets ${\mathcal P}_0$ and ${\mathcal Q}_0$.
In fact, to any solution  $Q=\left[(I-\Pi) Q_0^{-1} (I-\Pi)\right]\pinv\in {\mathcal Q}_0$ there corresponds a solution $P=\left[\Pi P_0^{-1} \Pi\right]\pinv\in {\mathcal P}_0$. To see this, just note that the orthogonal projector $\Pi_{\script Y}$ onto ${\script Y}:={\script  X}^\perp$ is equal to $(I-\Pi)$, with $\Pi$ being the  orthogonal projector   onto ${\script  X} $. In this case we shall call $P$ and $Q$ {\em complementary solutions} of the LMI's (\ref{garep}) and  (\ref{gareq}). Indeed for complementary solutions we have
$$ \rank P + \rank Q =n \,.$$
Of course, when $A$ is not unmixed, a similar correspondence holds for any pair of families
 ${\mathcal P}_{\Delta}$ and ${\mathcal Q}_{\Delta'}$ of solutions of (\ref{garep}) and  (\ref{gareq}) respectively, where ${\mathcal P}_{\Delta}$ is the family corresponding to a certain $P_{\Delta}\in{\mathbb P}$  and
${\mathcal Q}_{\Delta'}$ is the family corresponding to   $Q_{\Delta'}:=P_{\Delta}^{-1}\in{\mathbb Q}$ (with $\Delta':=P_{\Delta}-P_0$).
\erem

\subsection{The case of $A$ non-singular: Riccati equations}

In  case of a non-singular $A$ matrix, equations (\ref{garep}) and (\ref{gareq}) reduce, respectively, to the following homogeneous algebraic Riccati equations (ARE)
\beq\label{RicP}
P = A P A\tp - A P C\tp(I+CPC\tp)^{-1} CPA\tp\\
\eeq
and 
\beq\label{RicQ}
Q = A\tp Q A - A\tp Q B (I+B\tp QB)^{-1} B\tp QA.\\
\eeq
The equivalence of the two representations is stated in  the following proposition.
\bprop\label{arevsgare} 
Let (\ref{realization-of-Q-2new}) be a minimal realization of a rational discrete-time all-pass function and assume that $A$ is non-singular. 
Then $P=P\tp$ is a solution of (\ref{garep}) if and only if it is a solution of (\ref{RicP}) and $Q=Q\tp$ is a solution of (\ref{gareq}) if and only if it is a solution of (\ref{RicQ}).
\eprop
\proof
We prove only the equivalence of (\ref{gareq})  and (\ref{RicQ}) as the other equivalence is dual.
Let $Q$ be a solution of (\ref{gareq}). 
Then  there exist $H\in \R^{m\times n}$ and $J\in\R^{m\times m}$  such that $N(Q)=[H \mid J]\tp[H\mid J]$.
In view of Theorem \ref{maintheorem-ap}, point \ref{point3}), $H(zI-A)^{-1}B+J$ is all-pass. After  eliminating  the non-observable part of this realization we obtain a minimal realization say $\bar{C}(zI-\bar{A})^{-1}\bar{B}+J$ of the same all-pass function where the $\bar{A}$ matrix clearly remains non-singular.
This, in particular implies that $J$ is also non-singular so that $I+B\tp Q B=J\tp J$ is strictly positive definite and hence invertible.
Then, ${\rm rank}[N(Q)]=m$ implies that the Schur complement of $I+B\tp Q B$ in $N(Q)$ vanishes which is equivalent to $Q$ being a solution of (\ref{RicQ}).\\
Conversely, let $Q=Q\tp$ be  an arbitrary solution  of (\ref{RicQ}). To show that $Q$ satisfies the LMI (\ref{gareq}) 
it is enough  to show that $I+B\tp Q B$ is positive semi-definite and, hence, positive defnite.
In fact, in this case we can use, in the opposite direction, the previous argument based on the Schur complement.\\
The Riccati equation (\ref{RicQ}) can be written as
\beq
QA\inv =A\tp Q-A\tp QB(I+B\tp Q B)^{-1}B\tp Q
\eeq
from which it is easy to see that $\ker(Q)$ is   $A\inv$-invariant and hence $A$-invariant.
Select a basis where $A$ has the form shown in the right-hand side of (\ref{a-inv-sp-formnew}), $B=\bmat B_1\\B_2\emat$ is partitioned conformably and
$Q$ has   the same structure of the right-hand side of (\ref{q01}) where $P_{11}$ is non singular so that  $Q_{11}:= P_{11}^{-1}$ is also  non-singular.
Then substituting $Q= \diag \{Q_{11},\; 0\}$ into (\ref{RicQ}) it is immediate to see that $P_{11}^{-1}$ satisfies 
\beq\label{pre-Ric-S1new}
P_{11}^{-1}=A_1\tp P_{11}^{-1} A_1  - A_1\tp P_{11}^{-1}B_1(I+B_1\tp P_{11}^{-1}B_1)^{-1}B_1\tp P_{11}^{-1} A_1
\eeq
so that, by using the Sherman-Morrison-Woodbury formula, we get
$
A_1 P_{11} A_1\tp=P_{11} +B_1B_1\tp.
$
Since $(A,B)$ is reachable, $(A_1,B_1)$ is also reachable and Theorem \ref{maintheorem-ap}, point \ref{point4}), 
implies that $I+B_1\tp Q_{11} B_1$ is positive semidefinite. Observing that $I+B\tp Q B=I+B_1\tp Q_{11} B_1$
concludes the proof.
\qed
Notice that the ARE's (\ref{RicP}) and (\ref{RicQ}) do not impose explicitly any positivity condition: the previous result shows that these conditions are automatically met when $A$ is non-singular.
On the contrary, when $A$ is singular, it seems that one needs to impose explicitly the positivity condition in (\ref{garep}) and (\ref{gareq}):  this may be merely due to a technical difficulty and we conjecture that
the LMI (\ref{garep}) has the same solution set of  the equation ${\rm rank}[M(P)]=m$ and dually, the LMI  (\ref{gareq}) has the same solution set of equation ${\rm rank}[N(Q)]=m$.\\ 
As a direct consequence of Theorem \ref{thmhgare} and Proposition \ref{arevsgare}, we have the following corollary.
\bcor\label{corolRic-vero}   Let 
(\ref{realization-of-Q-2new}) be a minimal realization of a  rational bi-proper discrete-time all-pass function. Then
\begin{enumerate}
\item
The unique solution $P_0=P_0\tp$ of (\ref{eqforP}) is also a non-singular solution of the homogeneous Riccati equation (\ref{RicP}).
This solution generates the family ${\mathcal P}_0$ of symmetric solutions of (\ref{RicP}) as described by equation (\ref{paramgenforsetutte}), where  $P_{\Delta}=P_0$ and where $\Pi $ is the orthogonal projector onto an $A\tp$-invariant subspace ${\script  Y}$.
The elements $P=P\tp$ of this family are  in a one-to-one correspondence with the set of $A\tp$-invariant subspaces.\\
If $A$ is unmixed then $P_0$ is the only non-singular solution of (\ref{RicP})
and ${\mathcal P}_0$ is the set of all  solutions of (\ref{RicP}).
 
\item
The unique solution $Q_0=Q_0\tp$ of (\ref{eqforQ}) is also a non-singular solution of the homogeneous Riccati equation  (\ref{RicQ}).
This solution generates the family ${\mathcal Q}_0$ of symmetric solutions of (\ref{RicQ}) as described by equation (\ref{paramgenforsetutteq}), where  $Q_{\Delta}=Q_0$ and where $\Pi $ is the orthogonal projector onto an $A $-invariant subspace ${\script  X}$.
The elements of this family are  in a one-to-one correspondence with the set of $A$-invariant subspaces  ${\script  X}$.\\
If $A$ is unmixed then $Q_0$ is the only non-singular solution of (\ref{RicQ})
and ${\mathcal Q}_0$ is the set of all   solutions of (\ref{RicQ}).
\end{enumerate}
\ecor

\section{Factorization of all-pass functions}

In this section we discuss a remarkable relation between solutions of the constrained LMI's  (or ARE) and all pass divisors. The background facts are established in  the following lemma.

\blem\label{lembicambiobase}
Let (\ref{realization-of-Q-2new}) be a minimal realization of a  rational discrete-time all-pass function and 
let  $Q_0$ be the unique  solution of (\ref{eqforQ}).
Let $P\in{\mathcal P}_0$ and let  $Q\in{\mathcal Q}_0$ be  the complementary solution of (\ref{gareq}) associated to $P$ 
in the sense  described in Remark \ref{rembijcpq}.
Let ${\script X}:= \ker Q$ be the $A$-invariant subspace corresponding to $Q$ and
 ${\script Y}:= \ker P ={\script X}^\perp$ be the $A\tp$-invariant subspace corresponding to $P$. 
Then, one can select a basis such that, ${\script X}, {\script Y}, A,B,C,Q,P$ and  $Q_0$ have the following structure
\beq\label{sts1t1nb}
{\script X}=\im\bmat 0\\I\emat, \ {\script Y}=\im\bmat I\\0\emat,\eeq
\beq\label{stabcnb}
A=\bmat A_r&0\\A_{21}&A_l\emat,\ B=\bmat B_1\\B_2\emat,\  C=[C_1\mid C_2],
\eeq
\beq\label{stq1p1qnb}
 Q =\bmat P_{11}^{-1}&0\\0&0\emat, \  P =\bmat 0&0\\0&Q_{22}^{-1}\emat,\  Q_0=\bmat P_{11}^{-1}&0\\0&Q_{22}\emat.
 \eeq
\elem
\proof
Perform a preliminary change of basis as in equation (\ref{s1newbasisnew}) of the proof of Theorem \ref{thmhgare} (but  now use a slightly different notation) so that ${\script X}, {\script Y}$ are given by (\ref{sts1t1nb}), and $A$ has the block-triangular structure
$A=\bmat A_r&0\\\bar{A}_{21}&A_l\emat$.
In this basis,   partition $Q_0$ and $P_0=Q_0^{-1}$ as
$Q_0=\bmat Q_{11}&Q_{12}\\Q_{12}\tp&Q_{22}\emat$ and   $P_0=\bmat P_{11}&P_{12}\\P_{12}\tp&P_{22} \emat$. We have already proved  that in this basis $Q =\bmat P_{11}^{-1}&0\\0 &0\emat$.
Considering  (\ref{paramgenforsetutte}), where we set $P_\Delta=Q_0^{-1}$, and $\Pi $ is the orthogonal projector onto 
${\script Y}$, we   see that in the same basis $P=\bmat 0&0\\0&Q_{22}^{-1}\emat$.
Partition $B$ and $C$ conformably as
$B=\bmat B_1\\\bar{B}_2\emat,$ and $C=[\bar{C}_1\mid C_2]$.
From the first of (\ref{eqforQ})  it   follows that 
\beq\label{lyapridottaperq}
A_l\tp Q_{22} A_l - Q_{22} =C_2\tp C_2.
\eeq
Since $(A,C)$ is observable, $(A_2,C_2)$ is observable  as well so that from
(\ref{lyapridottaperq}) it follows that $Q_{22}$ is non-singular, \cite[Lemma 3.1]{Ferrante-Ntog-Automatica-13}.
Since $Q$ and $Q_{22}$ are non-singular the Schur complement $Q_{11}-Q_{12}Q_{22}\inv Q_{12}\tp$ is also non-singular and $P_{11}=(Q_{11}-Q_{12}Q_{22}\inv Q_{12}\tp)^{-1}$.
Perform now  a further change of basis induced by  $T=\bmat I&0\\-Q_{22}^{-1}Q_{12}\tp&I\emat$.
Although the two subspaces (\ref{sts1t1nb}) are no longer orthogonal they are still in direct sum, matrix $A$ in (\ref{stabcnb}) is modified just by changing $\bar A_{21}$ into 
$A_{21}:= \bar{A}_{21}+Q_{22}\inv Q_{12}\tp A_1 - A_2 Q_{22}\inv Q_{12}\tp$, and   in (\ref{stabcnb}) we have
$B_2:=\bar{B}_2+Q_{22}\inv Q_{12}\tp B_1$, and  $C_1:=\bar{C}_1-C_2 Q_{22}\inv Q_{12}\tp$.
\qed

\bthm \label{theofacorizationallpass}
Let (\ref{realization-of-Q-2new})
be a minimal realization of a rational discrete-time all-pass function. 
Let ${\mathcal P}_0$ be the family of solutions of (\ref{garep})  associated to the (unique) solution $P_0$ of (\ref{eqforP}) and  ${\mathcal Q}_0$ be the family of solutions of 
(\ref{gareq}) associated to the (unique)  solution $Q_0$ or (\ref{eqforQ}), as described in Remark \ref{remdefcalpq}.\footnote{As already observed, under the additional assumption that $A$ is unmixed, ${\mathcal P}_0$ is the family of {\em all} symmetric solutions of (\ref{garep}) and ${\mathcal Q}_0$ is the family of {\em all} symmetric solutions of (\ref{gareq}).} 
\begin{enumerate}
\item
For each $P\in{\mathcal P}_0$, let $G$ and $L$ be such that $[G\tp\mid L\tp]$ has full row-rank and
\beq\label{factormp1}
M(P)=[G\tp\mid L\tp]\tp [G\tp\mid L\tp].
\eeq
Then  
\beq\label{realization-of-Ql}
 Q_{L}(z):= C(zI-A)^{-1}G +L
\eeq
is a (non-minimal) realization of a left all-pass divisor of $Q(z)$.
The McMillan degree $n_l$ of  $Q_{L}(z)$ is equal to the rank of $P$.

Conversely, any left all-pass divisor of $Q(z)$ is given by  (\ref{realization-of-Ql}), 
where $[G\tp\mid L\tp]$ has full row-rank and satisfies (\ref{factormp1}) for a suitable $P\in{\mathcal P}_0$.
\item
For each $Q\in{\mathcal Q}_0$, let $H$ and $J$ be such that $[H\mid J]$ has full row-rank and
\beq\label{factornq1}
N(Q)=[H\mid J]\tp [H\mid J].
\eeq
Then  
\beq\label{realization-of-Qr}
 Q_{R}(z):= H(zI-A)^{-1}B +J
\eeq
is a (non-minimal) realization of a right all-pass divisor of $Q(z)$.
The McMillan degree $n_r$ of  $Q_{R}(z)$ is equal to the rank of $Q$.

Conversely, any right all-pass divisor of $Q(z)$ is given by  (\ref{realization-of-Qr}), 
where $[H\mid J]$ has full row-rank and satisfies (\ref{factornq1}) for a suitable $Q\in{\mathcal Q}_0$.
\end{enumerate}

\ethm

\proof
We prove only point 1) as point 2) is dual.
We first observe that $Q_{L}(z)$ is all-pass; in fact, $P$ is a solution of (\ref{eqforP}) associated to the quadruple $(A,G,C,L)$ so that in view of point \ref{point3}) of Theorem \ref{maintheorem-ap}, $Q_{L}(z)$ is all-pass.

Let $P\in{\mathcal P}_0$ and   $Q\in {\mathcal Q}_0$ be complementary solutions of (\ref{garep}) and (\ref{gareq}), respectively, as described in Remark \ref{rembijcpq}.
Select a basis as in Lemma \ref{lembicambiobase} so that ${\script X}, {\script Y}, A,B,C,Q,P$ and  $Q_0$ have the  structure described in (\ref{sts1t1nb}), (\ref{stabcnb}) and (\ref{stq1p1qnb}).
In the chosen basis,   compute $M(P)$ to obtain
\beq
M(P)=\bmat 0 & 0 & 0\\
0& A_l Q_{22}^{-1} A_l\tp -Q_{22}^{-1} &  A_l Q_{22}^{-1} C_2\tp\\
0 & C_2 Q_{22}^{-1} A_l\tp & C_2 Q_{22}^{-1} C_2\tp +I\emat
\eeq
so that $G$ must have the block structure,
$G=\bmat 0\\G_{l}\emat$ and $Q_{L}(z)$ defined in (\ref{realization-of-Ql}) has the following realization
\beq\label{minrealql1}
Q_{L}(z)=C_2(zI-A_l)^{-1}G_{l} +L
\eeq
(it could be shown that this realization is minimal but we will find that this result comes as a byproduct 
at the end of the proof). 
Observe now that $(A,C)$ is observable so that $(A_l,C_2)$ is observable as well.
Now since $Q_{22}^{-1}$ is a solution of (\ref{eqforP}) associated with the quadruple $(A_l,G_l,C_2,L)$, then  $Q_{22}$ must be  a solution of (\ref{eqforQ}) associated with the same quadruple.
In particular, from the second of equations (\ref{eqforQ}) we get
$A_l\tp Q_{22} G_l = C_2\tp L$ which may be rewritten as
$[A_l\tp Q_{22} \mid -C_2]\bmat G_l \\ L\emat=0$.
In this factorization, the matrix $[A_l\tp Q_{22} \mid -C_2]$ has full row-rank; in fact,
 $(A_l,C_2)$ is observable so that $[A_l  \mid C_2]$ has full row-rank; hence 
$[A_l  \mid -C_2]$ has full row-rank as well; furthermore since $Q_{22}$ is non-singular  also $[A_l\tp Q_{22} \mid -C_2]$ has full row-rank.
The right matrix $\bmat G_l \\ L\emat$ has full column-rank; in fact, we have already observed that its transpose has full row-rank.
In conclusion, $[A_l\tp Q_{22} \mid -C_2]\in\R^{n\times (n+m)}$ has rank $n$ so that its kernel has dimension $m$ and hence the $m$ linearly independent columns of the matrix $\bmat G_l \\ L\emat$ are a basis for $\ker[A_l\tp Q_{22} \mid -C_2]$.

Now use the fact that $Q_0$ is a solution of equations (\ref{eqforQ}) associated with the quadruple $(A,B,C,D)$. From the lower block of the second of these equations, we get
\beq
[A_l\tp Q_{22} \mid -C_2]\bmat B_2 \\ D\emat=0.
\eeq
Similarly, from the left-lower block of the first of the same equations, we get
\beq
[A_l\tp Q_{22} \mid -C_2]\bmat A_{21} \\ C_1\emat=0.
\eeq
Hence, there exist matrices $D_r$ and $C_r$ such that
\beq\label{67eqpazzesche}
\bmat B_2 \\ D\emat=\bmat G_l \\ L\emat D_r;\qquad 
\bmat A_{21} \\ C_1\emat=\bmat G_l \\ L\emat C_r.
\eeq
It is now a matter of direct computation to see that 
\bea\nn
Q(z)&\!\!\!\!=\!\!\!\!&[L C_r\mid C_2]\left(zI-\bmat A_r&0\\G_l C_r&A_l\emat\right)^{-1}\bmat B_1 \\ G_lD_r \emat C_r\\
\nn
&& + L  D_r\\
\nn
&\!\!\!\!=\!\!\!\!&[C_2(zI-A_l)^{-1}G_{l} +L][C_r(zI-A_r)^{-1}B_{1} +D_r]\\
&\!\!\!\!=\!\!\!\!&Q_{L}(z) \hat{Q}_R(z)
\label{factorizationl1r}
\eea
where we have introduced  the rational function $\hat{Q}_R(z):=C_r(zI-A_r)^{-1}B_{1} +D_r$.
Note that, since $Q(z)$ and $Q_{L}(z)$ are all-pass, $\hat{Q}_R(z)$ is necessarily all pass.
To show that  $Q_{L}(z)$ is a left divisor of $Q(z)$ it remains only to observe that
$Q(z)$ has a minimal realization of dimension $n$ and that $n=n_l+n_r$ where $n_l$ is the dimension of $A_l$ and $n_r$ is the dimension of $A_r$. As a byproduct,  (\ref{minrealql1}) is a minimal realization of $Q_{L}(z)$  and $C_r(zI-A_r)^{-1}B_{1} +D_r$ is a minimal realization of $\hat{Q}_R(z)$.
Finally, by construction, the McMillan degree $n_l$ of $Q_{L}(z)$   equals  the dimension of $Q_{22}$ or, equivalently,   the rank of $P$.

Conversely, let $Q(z)=\hat{Q}_L(z)\hat{Q}_R(z)$   with  $\hat{Q}_L(z):= C_l(zI-A_l)^{-1}B_l+D_l$, and $\hat{Q}_R(z):= C_r(zI-A_r)^{-1}B_r+D_r$, being  minimal realizations of   all-pass functions   and assume that the McMillan degree of $Q(z)$ equals the sum of the McMillan degrees of $\hat{Q}_L(z)$ and $\hat{Q}_R(z)$. 
  Then, up to a change of basis which does not affect the result that we need to establish, we have that
\beq\label{abseries}
A=\bmat
A_r & 0\\ B_lC_r & A_l\emat;\quad B=\bmat B_r\\B_lD_r\emat.
\eeq
\beq\label{cdseries}
C=[D_lC_r\mid C_l];\quad D=D_lD_r.
\eeq
Hence, without loss of generality, we assume that the matrices $A,B,C,D$ of (\ref{realization-of-Q-2new}) have the expressions (\ref{abseries}) and (\ref{cdseries}).
Since $\hat{Q}_L(z)$ and $\hat{Q}_R(z)$ are all-pass functions, there exist an invertible matrix
$P_l$  solving equations (\ref{eqforP}) associated with the quadruple $(A_l,B_l,C_l,D_l)$ and
and an invertible matrix $P_r$  solving equations (\ref{eqforP}) associated with the quadruple $(A_r,B_r,C_r,D_r)$.
By exploiting \eqref{abseries} and \eqref{cdseries}, it is straightforward to check that, in the selected basis, $\diag(P_r,P_l)$ is the (unique) solution of (\ref{eqforP}) associated with  the quadruple $(A,B,C,D)$. Hence, we have
$
P_0=\bmat P_r&0\\0&P_l\emat.
$
Let ${\script  Y}=\im\bmat I\\0\emat$ be an $A\tp$-invariant subspace so that
\beq\label{pcostruitaapposta}
P:=\left[(I-\Pi) P_0^{-1} (I-\Pi)\right]\pinv=\bmat 0&0\\0&P_l\emat\in{\mathcal P}_0.
\eeq
By direct computation, it is also straightforward to check that 
\beq\label{factcostruitaapposta}
M(P)=[G\tp\mid L\tp]\tp [G\tp\mid L\tp]
\eeq
where $G:= \bmat 0\\B_l\emat$ and $L:= D_l$.
Now     define, as in (\ref{realization-of-Ql}), the left factor $Q_{L}(z)$ associated with $P$, $G$ and $L$ given by (\ref{pcostruitaapposta}) and (\ref{factcostruitaapposta}). By eliminating non-reachable part of this $Q_{L}(z)$, we see  that
$Q_{L}(z)=\hat{Q}_{L}(z)$. 
\qed

%%come dÕaccordo, ecco il pezzo da aggiungere di peso dopo la dimostrazione del Theorem 4.1.

Theorem \ref{theofacorizationallpass} provides a one to one correspondence between the family ${\mathcal P}_0$ of solutions of (\ref{garep}) and left all-pass factors of $Q(z)$ defined up to multiplication from the right side by a constant orthogonal matrix $U$.
Similarly, Theorem \ref{theofacorizationallpass} also provides a one to one correspondence between the family ${\mathcal Q}_0$ of solutions of (\ref{gareq}) and right factors of $Q(z)$ defined up to multiplication from the left side by a constant orthogonal matrix $U$.
On the other hand, a left factor $Q_{L}(z)$ of $Q(z)$ is associated with a right factor $Q_{R}(z)$ by the factorization relation $Q(z)=Q_{L}(z) Q_{R}(z)$. Given a factoriaztion of this type,
it is natural to ask what is the relation between the solution $P\in{\mathcal P}_0$ associated with $Q_{L}(z)$ and the solution $Q\in{\mathcal Q}_0$ associated with the corresponding $Q_{R}(z)$. The following result addresses this question and shows that $P$ and $Q$   are related by the same bijective correspondence introduced in Remark \ref{rembijcpq}.

\bprop
Let (\ref{realization-of-Q-2new})
be a minimal realization of a rational discrete-time all-pass function and let
$Q(z)=Q_{L}(z) Q_{R}(z)$ be a minimal factorization of $Q(z)$.
The matrices $P\in{\mathcal P}_0$ and $Q\in{\mathcal Q}_0$ associated with $Q_{L}(z)$ and $Q_{R}(z)$, respectively, by Theorem \ref{theofacorizationallpass} satisfy the relation $\ker [P]=(\ker[Q])^\perp$ and are therefore a complementary pair.
\eprop
\proof
As in the proof of theorem \ref{theofacorizationallpass}, let $P\in{\mathcal P}_0$ and let $Q\in {\mathcal Q}_0$ be the corresponding solution of (\ref{gareq}) as described in Remark \ref{rembijcpq}, i.e. the only element of ${\mathcal Q}_0$ such that
$\ker [P]=(\ker[Q])^\perp$.
We select a basis as in Lemma \ref{lembicambiobase} so that ${\script  X}, {\script Y}, A,B,C,Q,P$ and  $Q_0$ have the  structure described in (\ref{sts1t1nb}), (\ref{stabcnb}) and (\ref{stq1p1qnb}).
  Consider a left factor $Q_L(z)$ associated with $P$: as we have seen 
in the proof of Theorem \ref{theofacorizationallpass},  the corresponding right factor 
$\hat{Q}_R(z)$ (that satisfies equation (\ref{factorizationl1r})) has a minimal realization of the form
$\hat{Q}_R(z)=C_r(zI-A_r)^{-1}B_{1} +D_r$.
Let $P_r$ be the solution of (\ref{eqforP}) associated with the quadruple $(A_r,B_1,C_r,D_r)$.
By taking (\ref{stabcnb}) and (\ref{67eqpazzesche}) into account, we easily see by a direct computation that 
$\diag(P_r,Q_{22}^{-1})$ is the solution of (\ref{eqforP}) associated to the quadruple $(A,B,C,D)$.
Since the solution $P=\diag(P_{11},Q_{22}^{-1})$ of this equation is unique, we have $P_r=P_{11}$. On the other hand, we know that the right factor $Q_{R}(z)$ associated with the matrix $Q$ is given by (\ref{realization-of-Qr}) and, by duality, has a minimal realization of the form $Q_{R}(z)=H_r(zI-A_r)^{-1}B_{1} +J$.
Now we compare the all-pass functions $\hat{Q}_R(z)$ and  $Q_R(z)$ and we see that they have the same state and input matrices and that the solutions of the  equation (\ref{eqforP}) associated to the minimal quadruple
$(A_r,B_1,C_r,D_r)$ and of the equation (\ref{eqforP}) associated to the minimal quadruple
$(A_r,B_1,H_r,J)$, coincide.
Hence, $Q_R(z)$ and  $\hat{Q}_R(z)$ differ for multiplication on the left side by a constant orthogonal matrix.  
\qed

In the case when $Q(z)$ is bi-proper | or, equivalently, $A$ and $D$ are non-singular   we know that 
(\ref{garep}) and  (\ref{gareq}) reduce to ARE's. Moreover, for any given solution $P$ of (\ref{garep}), (or, equivalently, of (\ref{RicP})) we can provide an explicit expression for the matrices $G$ and $L$ by solving (\ref{factormp1}).
The following corollary connects solutions of ARE's and all-pass factorizations. 

\bcor \label{theoDivallpassbiproprie}   
Let (\ref{realization-of-Q-2new})
be a minimal realization of a rational bi-proper discrete-time all-pass function. 
Let ${\mathcal P}_0$ be the family of solutions of (\ref{RicP})  associated with the solution $P$ of (\ref{eqforP}) and  ${\mathcal Q}_0$ be the family of solutions of 
(\ref{RicQ}) associated with the solution $Q$ or (\ref{eqforQ}), as described in Corollary \ref{corolRic-vero}.\footnote{Similarly to the general case, under the additional assumption that $A$ is unmixed, ${\mathcal P}_0$ is the family of {\em all} symmetric solutions of (\ref{RicP}) and ${\mathcal Q}_0$ is the family of {\em all} symmetric solutions of (\ref{RicQ}).}

\begin{enumerate}
\item
For each $P\in{\mathcal P}_0$, the function
\beq\label{realization-of-Ql-bipropr}
 Q_{L}(z):= C(zI-A)^{-1}G+L,
\eeq
with
\beq\label{eqforDB}
\left\{
\begin{array}{l}
L:=(I+CP C\tp)^{1/2}\\
G:=A P C\tp L \mtp
\end{array}
\right.
\eeq
is a (non-minimal) realization of a left all-pass divisor of $Q(z)$.\\
Conversely, any left all-pass divisor of $Q(z)$ is given   up to multiplication from the right side by a constant orthogonal matrix   by  (\ref{realization-of-Ql-bipropr}), (\ref{eqforDB}).

\item
For each $Q\in{\mathcal Q}_0$, the function
\beq\label{realization-of-Qr-bipropr}
 Q_{R}(z):= H(zI-A)^{-1}B+J,
\eeq
with
\beq\label{eqforDC}
\left\{
\begin{array}{l}
J :=(I+B\tp Q B)^{1/2}\\
H:=J  \mtp B\tp Q A
\end{array}
\right.
\eeq
is a (non-minimal) realization of a right all-pass divisor of $Q(z)$.\\
Conversely, any right all-pass divisor of $Q(z)$ is given | up to multiplication on the left side by a constant orthogonal matrix   by  (\ref{realization-of-Qr-bipropr}), (\ref{eqforDC}).
\end{enumerate}
\ecor

\section*{Summary  and possible generalizations }
In this paper we have provided a completely general  characterization   of  discrete-time all-pass matrix  functions in the same spirit  of the continuous-time      result of   Glover's  \cite[Theorem 5.1 ]{Glover-84}. Applications to some class of  LMI's, to homogeneous Riccati equations and to the factorization of all-pass functions are discussed. The characterization  is  presented for square    all-pass matrix  functions but a generalization to non-square functions can be pursued along the same lines.

\appendices
\section{Factorization of all-pass functions which are singular at infinity}
\blem\label{primolemmafactoriallp}
Let $Q(z)$
be an $m\times m$ rational proper discrete-time all-pass function. 
Then $Q(z)$ can be written as 
\beq\label{facaps+ns}
Q(z)=Q_0(z)\bar{Q}_1(z)\bar{Q}_2(z)\dots \bar{Q}_k(z)
\eeq
where $Q_0(z)$ is a rational discrete-time all-pass function such that
$Q_0(\infty)$ is non-singular and   the $\bar{Q}_i(z)$'s are rational proper all-pass functions (whose only pole is in the origin) having a  realization of the following form
\beq\label{fdeiqbari}
\bar{Q}_i(z)=\bmat I_{m-p_i}&0\\0 & 0\emat U_i +\bmat 0\\I_{p_i}\emat (zI_{p_i}-0)^{-1}[0\mid I_{p_i}]U_i
\eeq
where $U_i$ is a constant orthogonal matrix.
\elem
\proof
Consider a minimal realization $Q(z)=C(zI-A)^{-1}B+D$.
If $D$ is non-singular, $Q_0(z)=Q(z)$ and we are done.
If $D$ is  singular, we resort to the  Silverman
algorithm as described in \cite{Ferrante-P-P-02-LAA}.
Assume the matrix $D$ has $q_1$ linearly independent columns, with $0\leq
q_1< m$. Let $V_1$ be an orthogonal matrix such that
$DV_1=\bmat D_{11} & 0\emat$, with $D_{11}\in\R^{m\times q_1}$ being
full column rank. Let us partition $BV_1=\bmat B_{11} & B_{12}\emat$
conformably, obtaining the following block structure,
\beq
\tilde{Q}_1 (z):=Q(z)V_1 = C(zI-A)^{-1}\bmat B_{01} &
B_{02}\emat+\bmat D_{01} & 0\emat,
\eeq
and let
\beq \hat Q_{1}(z):=\tilde{Q}_1 (z) \bmat I_{q_1}&0\\0&zI_{m-q_1}\emat.
\eeq
Clearly, $\hat Q_{1}(z)$ is all-pass as it is the product of all-pass functions.
Moreover, $\hat Q_{1}(z)$ can be written as
$$
\hat Q_{1}(z)=[\hat Q_{11}(z)\mid \hat Q_{12}(z)]
$$
where
$$\hat Q_{11}(z):= D_{11} + CB_{11} z^{-1} + CAB_{11} z^{-2} +
\ldots 
$$
and
$$
\hat Q_{12}(z):=CB_{12} + CAB_{12} z^{-1} + CA^2B_{12} z^{-2}+ \ldots
$$
so that $\hat Q_{1}(z)$ has the following realization
\beq
\hat Q_{1}(z)=  C(zI-A)^{-1}\bmat B_{11} & AB_{12}\emat+\bmat D_{11} &
CB_{12}\emat.
\eeq

At this point, either $\bmat D_{11} & CB_{12}\emat$ is
right-invertible, or we may iterate the above procedure by introducing
another orthogonal matrix $V_2$, such that
$$
\bmat D_{01} & CB_{02}\emat V_2=\bmat D_{21} & 0\emat,
$$
with $D_{21}\in\R^{m\times q_2}$ of full column rank and $q_2\geq q_1$; we
 define the new all-pass function
\beq
\tilde{Q}_2 (z):=\hat{Q}_1(z)V_2=C(zI-A)^{-1}\bmat  B_{21} & B_{22}\emat +
\bmat  D_{21} & 0\emat,
\eeq
where $\bmat  B_{21} & B_{22}\emat=\bmat  B_{11} & AB_{12}\emat
V_2$.

Since $Q(z)$ is all-pass, it has full rank (as a rational matrix function) and hence, after a finite number of steps (say $k$) of the above
procedure, we get a rational proper all pass function
\beq \label{Qk(z)}
\tilde{Q}_k (z)=  Q(z) \prod_{i=1}^{k}V_i\bmat I_{q_i}&0\\0&zI_{m-q_i}\emat,
\eeq
such that $\tilde{Q}_k (\infty)$ is non-singular.
Now we set
$Q_0(z):=\tilde{Q}_k (z)$, so that
\beq \label{Q(z)}
 Q(z) =  Q_0(z) \left[\prod_{i=1}^{k}V_i\bmat I_{q_i}&0\\0&zI_{m-q_i}\emat\right]\inv.
\eeq
Finally, by setting
 $p_i:= q_{k+1-i}$, $i=1,2,\dots,k$, and $U_i:=V_{k+1-i}\tp$, $i=1,2,\dots,k$, and observing that 
 \beq
\bmat I_{p_i}&0\\0&zI_{m-p_i}\emat\inv =\bmat I_{m-p_i}&0\\0 & 0\emat  +\bmat 0\\I_{p_i}\emat (zI_{p_i}-0)^{-1}[0\mid I_{p_i}]
\eeq
we obtain (\ref{facaps+ns}) and (\ref{fdeiqbari}). 
\qed

\blem\label{secondolemfacapf}
Let $Q(z)$
be an $m\times m$ rational proper discrete-time all-pass function factored as in (\ref{facaps+ns}).
Consider a reachable realization 
\beq
Q_{i}(z)=C_{i}(zI-A_{i})B_i +D_i
\eeq
of $Q_{i}(z):= Q_0(z)\bar{Q}_1(z)\bar{Q}_2(z)\dots \bar{Q}_{i}(z)$. Partition $B_i$ and $D_i$ as $B_i=[B_{i,1}\mid B_{i,2}]$ and $D_i=[D_{i,1}\mid D_{i,2}]$, where $B_{i,1}$ and $D_{i,1}$ have $m-p_{i+1}$ columns.
Then a reachable realization  of $Q_{i+1}(z):=Q_{i}(z)\bar{Q}_{i+1}(z)$ is given by
\bea\nn
Q_{i+1}(z)\!\!\!&\!\!\!\!\!\!\!=\!\!\!\!\!\!\!&\!\!\![D_{i,2}\mid C_i]\left(zI-\bmat 0&0\\B_{i,2}&A_i\emat\right)^{-1}\bmat 0 & I\\B_{i,1}&0\emat U_{i+1}\\
\label{realqi+1}
&& +[D_{i,1}\mid 0]U_{i+1}.
\eea
\elem
\proof
The realization (\ref{realqi+1}) is the result of a direct computation.
The fact that this realization is reachable may be easily seen  by using the PBH test.
In fact, as a consequence of the fact that $[A_i-\lambda I \mid B_{i,1}\mid B_{i,2}]$
has full row-rank for all $\lambda\in\C$, we immediately see that also
$$
\bmat -\lambda I&0&0 & I\\B_{i,2}&A_i-\lambda I&B_{i,1}&0\emat  $$
has full row-rank for all $\lambda\in\C$.
\qed

\section{Proof of symmetry of $T$}\label{B}
Somehow in the same spirit of \cite{Glover-84}, 
 we shall show that
 \beq
 U:=T^{-1} T\tp
 \eeq
 satisfies
 \bea
 \stoa\label{simia}
&& A= U^{-1} A U,\\
 \stob\label{simib}
 && B = U^{-1} B,\\
 \stoc\label{simic}
&& C = C U.
 \eea\rst
This means that $U$ is a similarity transform that leaves unchanged the triple $(A,B,C)$ of the system.
Since $(A,B,C)$ is, by assumption, a minimal realization, this means that $U=I$, or that
$T=T\tp$.

We start with (\ref{simic}).
Solving (\ref{eqiuvb}) and (\ref{eqiuvc}) for $B$
we get
\beq\label{exb1}
B=T^{-1}A^{-\top} C\tp D
\eeq
and
\beq\label{exb2}
B=AT^{-\top} C\tp D^{-\top}
\eeq
By inserting in the latter the expression of $D^{-\top}$ obtained by transposing (\ref{eqiuvd}), we get 
\beq\label{exb3}
B=AT^{-\top}  C\tp D - A T^{-\top}C\tp CA^{-1} B
\eeq
Now we take the inverse of both sides of (\ref{eqiuva}) and use the Sherman-Morrison-Woodbury formula thus obtaining
\bea\nn
T^{-1}A\tp T &=&A^{-1}+A^{-1} B(D-CA^{-1} B)^{-1} CA^{-1}\\
\label{SMWFper}
&=&A^{-1}+A^{-1} BD\tp CA^{-1}.
\eea
From (\ref{exb2}) we get $BD\tp=AT^{-\top} C\tp $ which, plugged into the left-hand side of (\ref{SMWFper}), yields
\beq
T^{-1}A\tp T=A^{-1}+T^{-\top}  C\tp CA^{-1}.
\eeq
The latter provides an expression for $T^{-\top} C\tp CA^{-1}$ which, plugged into the left-hand side of (\ref{exb3}) gives
\beq\label{exb4}
B=AT^{-\top}  C\tp D +B - A T^{-1}A\tp TB 
\eeq
so that
$
T^{-\top}  C\tp D=T^{-1}A\tp TB,
$
or
$
B=T^{-1}A^{-\top}TT^{-\top}C\tp D.
$
By comparing the latter with (\ref{exb1}), we eventually get
\beq
C\tp = TT^{-\top}C\tp
\eeq
which, by  recalling that 
$U:=T^{-1} T\tp,$ readily implies (\ref{simic}).

We now use a dual argument to obtain (\ref{simib}).
Solving (\ref{eqiuvb}) and (\ref{eqiuvc}) for $C$
we get
\beq\label{exc1}
C= D B\tp A\mtp T
\eeq
and
\beq\label{exc2}
C=D\mtp B\tp T\tp A
\eeq
By inserting in the latter the expression of $D^{-\top}$ obtained by transposing (\ref{eqiuvd}), we get 
\beq\label{exc3}
C=D  B\tp T\tp A -CA\inv BB\tp T\tp A
\eeq
From (\ref{exc2}) we get $D\tp C=B\tp T\tp A $ which, plugged into the left-hand side of (\ref{SMWFper}), yields
$
T^{-1}A\tp T=A^{-1}+ A^{-1} BB\tp T\tp.
$
The latter provides an expression for $A^{-1} BB\tp T\tp$ which, plugged into the left-hand side of (\ref{exc3}) gives
\beq\label{exc4}
C=D  B\tp T\tp A +C -CT\inv A\tp T A
\eeq
so that
$
D  B\tp T\tp=CT\inv A\tp T,
$
or
$
C=D  B\tp T\tp T\inv A\mtp T.
$
By comparing the latter with (\ref{exc1}), we eventually get
\beq
B\tp = B\tp T\tp T\inv
\eeq
which, by  recalling that 
$U:=T^{-1} T\tp,$ readily implies (\ref{simib}).

We now prove   (\ref{simia}).
We multiply equation (\ref{eqiuva})
on the left side by $U^{-1}$ and on the right side by $U$. 
By taking into account (\ref{simib}) and (\ref{simic}), we get
\beq\label{afintocb}
U^{-1} A U = T\mtp A\mtp T\tp + BD\inv C.
\eeq
On the other hand, by transposing the first and the last member of (\ref{SMWFper}) and multiplying on the left side by $T\mtp$ and on the right side by $T\tp$  we get
\bea\nn
A&=&T\mtp A\mtp T\tp+T\mtp A\mtp C\tp D B\tp A\mtp T\tp\\
\nn
&=&T\mtp A\mtp T\tp+T\mtp A\mtp C\tp DD\inv D B\tp A\mtp T\tp.\\
\eea
Moreover, by inserting in the right-hand side of the latter the expressions of
$A\mtp C\tp D$ and  $D B\tp A\mtp$ obtained from (\ref{exb1}) and (\ref{exc1}), respectively, we get
\bea\nn
A &=& T\mtp A\mtp T\tp+ \underbrace{T\mtp T}_{U^{-1}}B D\inv C\underbrace{T^{-1} T\tp}_{U}\\
&=& T\mtp A\mtp T\tp+B D\inv C.
\eea
Finally, by comparing the latter with (\ref{afintocb}), we get (\ref{simia}).
\qed

\end{document}